\documentclass[notitlepage,a4paper]{article}

\usepackage[english,french]{babel}
\usepackage[T1]{fontenc}
\usepackage[dvips]{graphicx}
\usepackage{amsmath,amsfonts,amssymb,amsthm,latexsym,graphicx}
\usepackage{natbib}
\usepackage{color}
\usepackage[applemac]{inputenc}
\usepackage{pstricks}
\usepackage{pst-node}

\newtheorem{theorem}{Theorem}
\newtheorem{corollary}{corollary}

\newtheorem{example}{Example}
\newtheorem{e-proposition}{Proposition} \theoremstyle{definition}
\newtheorem{e-definition}{Definition} \theoremstyle{definition}
\newtheorem{remark}{Remark} 

\begin{document} 
\title{Orbit measures and interlaced determinantal point processes}
\author{M. Defosseux\\Laboratoire de Probabilit\'es et Mod\`eles Al\'eatoires}
\maketitle

\selectlanguage{english}

\begin{abstract}
We study random interlaced configurations on $\mathbb{N}\times \mathbb{R}$ considering the eigenvalues of the  main minors of Hermitian random matrices of the classical complex Lie algebras. We show that these random configurations are determinantal and give their correlation kernels.

\vskip 0.5\baselineskip

\end{abstract}

\section{Introduction}
This note announces the results of  \cite{Defosseux}. Baryshnikov \cite{Baryshnikov} has studied  the law of the eigenvalues of the main minors of a random matrix from the GUE, developing the connexions with uniform measure on Gelfand Cetlin patterns. In this paper we establish such connexions in a more general context, using a version  of Heckman's theorem \cite{Heckman}.  From this theorem and from classical branching rules, we deduce the joint law of the eigenvalues of the main minors of random Hermitian matrices of classical complex Lie algebras, with invariant law under some unitary transformations. This alows us to show that the associated point processes are determinantal. For the GUE minors process, this has been proved  by Johansson and Noordenstam \cite{Johansson} and  Okounkov and Reshetikhin \cite{Okounkov}  (see also,  very recently, Forrester and Nagao \cite{Forrester}). Thus, we generalise some of their results to all the  classical complex Lie algebras.
\section{Approximation of orbit measures}
\label{Approximation}
Let $G$ be a connected compact Lie group with Lie algebra $\mathfrak{g}$ and complexified Lie algebra $\mathfrak{g}_\mathbb{C}$. We choose a maximal torus $T$ of $G$ and we denote  by $\mathfrak{t}$ its Lie algebra. We equip $\mathfrak{g}$ with an Ad($G$)-invariant inner product $\langle.,.\rangle$ which induces a linear isomorphism between $\mathfrak{g}$ and its dual $\mathfrak{g}^*$ and intertwines the adjoint and the coadjoint action of $G$. We consider  the roots system $R$, i.e. all $\alpha\in \mathfrak{t}^*$ such that there is a non zero $X\in \mathfrak{g}_\mathbb{C}$ such that for all $H\in \mathfrak{t}$, $[H,X]=i\alpha (H)X$. We choose  the set $\Sigma$ of simple roots of $R$, $C(G)=\{\lambda\in \mathfrak{t}^*: \langle\lambda,\alpha \rangle >0 \textrm{ for all } \alpha\in \Sigma\}$ the corresponding Weyl chamber and $P^+(G)=\{\lambda\in \mathfrak{t}^*: 2\frac{\langle\lambda,\alpha \rangle}{\langle\alpha,\alpha \rangle}\in \mathbb{N}, \textrm{ for all } \alpha\in \Sigma_G\}$ the corresponding set of integral dominant weights. 

We consider  a connected compact subgroup $H$ of $G$ with Lie algebra $\mathfrak{H}$. We can choose a maximal torus $S$ of $H$ contained in $T$ and a corresponding set of integral dominant weights denoted by $P^+(H)$. For $x\in \mathfrak{g}^*$, let $\pi_H(x)$ be the orthogonal projection of $x$ on $\mathfrak{H}^*$.

For $\lambda\in P^{+}(G)$, we denote by $V_\lambda$ the irreducible  $\mathfrak{g}$-module with highest weight $\lambda$ and $\dim_G(\lambda)$ the dimension of $V_\lambda$.  For $\beta\in P^+(H)$ we denote by $m^\lambda_H(\beta)$ the multiplicity of the $\mathfrak{H}$-module with highest weight $\beta$ in the decomposition into irreducible components of $V_{\lambda}$ considered as an $\mathfrak{H}$-module. Rules giving the value of the multiplicities $m^\lambda_H$ are called \textit{branching rules}.

The intersection between the orbit of an element $x\in \mathfrak{g}^*$ under the coadjoint action of $G$ and the closure $\bar{C}(G)$ contains a single point that we call the radial part of $x$ and denote by $r_G(x)$. The same holds for $H$ and, for $x\in \mathfrak{g}^*$, we write $r_H(x)$ instead of $r_H(\pi_H(x))$.

The following proposition is a version of Heckman's theorem  \cite{Heckman} on asymptotic behaviour of multiplicities. As usal $\delta_x$ is the Dirac measure at $x$. \\
\begin{e-proposition} \label{theogeneral}  Let $x\in \bar{C}(G)$. 
Let $(\lambda_{n})_{n\ge 0}$ be a sequence of elements in
$P^{+}(G)$ and $(\epsilon_{n})_{n\ge 0}$ a sequence of positive real numbers such that $\epsilon_{n}$ converges to zero and 
$\epsilon_{n}\lambda_{n}$  converges to $x$, as $n$ tends to $+\infty$. Then   the sequence $(\mu_{n})_{n\ge 0}$ of probability measures on $\bar{C}(H)$ defined by 
$$\mu_{n}=\sum_{\beta\in P^+(H)}\frac{\dim_{H}(\beta)}{\dim_{G}(\lambda_{n})}m^{\lambda_{n}}_H(\beta)\delta_{\epsilon_n\beta}.$$
 converges to the law of
$r_{H}(Ad(g)x)$, with $g$ distributed according to the normalised Haar measure on $G$.
\end{e-proposition}
We deduce from this proposition and from classical branching rules the descriptions of convolutions or projections of $Ad(G)$-invariant measures on classical Lie algebras. Let us give a first application when $G$ is the symplectic unitary group $Sp(n)$. \\

\begin{example}. \label{sumC} This example concerns random variables which play the same role for $Sp(n)$ as the Laguerre ensemble for $U(n)$.
Let $(X_i)_{i\ge1},(Y_i)_{i\ge 1}$ be two independent sequences of random variables in $\mathbb{C}^n$  with independent  standard complex  Gaussian components. We consider $M_k=\sum_{i=1}^{k}S(X_i,Y_i)R(X_i,Y_i)^*$, $k\ge 1$, where, for $x=(x_1,\cdots,x_n),y=(y_1,\cdots,y_n)\in\mathbb{C}^n$,
\begin{displaymath}
S(x,y)=\left(%
\begin{array}{c}
s(x_1,y_1)\\
\vdots \\
s(x_n,y_n)
\end{array}%
\right), R(x,y)=\left(%
\begin{array}{c}
r(x_1,y_1)\\
\vdots \\
r(x_n,y_n)
\end{array}%
\right), s(a,b)=\left(%
\begin{array}{cc}
a & -b\\
\bar {b} & \bar{a}
\end{array}%
\right), r(a,b)=\left(%
\begin{array}{cc}
a & b\\
\bar {b} & -\bar{a}
\end{array}%
\right), a,b\in \mathbb{C}.
\end{displaymath} Let $\Lambda^{(k)}$ be the  positive eigenvalues of $M_k$. 
The process $(\Lambda^{(k)})_{k\ge 1}$ is an inhomogeneous Markovian process. Its transition kernel is deduced from Proposition \ref{theogeneral} considering tensor product by irreducible representations with highest weight proportional to the highest weight of the standard representation (see \cite{Defosseux}). \end{example} 
 
 \section{Eigenvalues of the main minors of classical Hermitian matrices}\label{Eigenvalues}
 
 \subsection{Classical Hermitian matrices} 
We study a  class of measures on generalised Gelfand Cetlin cones (see, e.g.,  Berenstein et Zelevinsky  \cite{Berenstein}). We denote by $M_n(\mathbb{C})$ (resp. $ M_n(\mathbb{R}))$ the set of $n\times n$ complex (resp. real) matrices and $I_n$  the unit matrix of  $M_n(\mathbb{C})$. We consider the classical compact groups. $G_a$ is the unitary group $U(n)=\{M\in M_{n}(\mathbb{C}): M^*M=I_n\}$, $G_b=SO(2n+1)$, $G_c$ is the unitary symplectic group written in a somewhat unusual way as $Sp(n)=\{M\in U(2n): M^{t}JM=J\}$ where all entries of $J=(J_{i,j})_{1\le i,j\le 2n}$ are zero except $J_{2i,2i-1}=-J_{2i-1,2i}=1$, $i=1,...,n$,  and $G_d=SO(2n)$, where $SO(k)$ is the orthogonal group $\{M\in M_{k}(\mathbb{R}): M^{*}M=I_{k}, \det(M)=1\}$. These compact groups correspond to the roots system of type $A,B,C,D$.  We write  $\mathfrak{g}_\nu$ their Lie algebras and $\mathcal{H}_\nu=i\mathfrak{g}_\nu$, $\nu =a,b,c,d$, the sets of associated  Hermitian matrices. 

The radial part of $M \in \mathcal{H}_\nu$ can be identified with the ordered eigenvalues of $M$ when $\nu=a$, resp. ordered positive eigenvalues of $M$ when $\nu=b,c$. 
When $\nu=d$, the Weyl chamber is $C=\{D(x): x\in \mathbb{R}^n, x_1 >... > x_{n-1} >\vert x_n\vert \}$ where $D(x)$ is the matrix in $ M_{2n}(\mathbb{C})$ all the entries of which are zero except $D(x)_{2k,2k-1}=-D(x)_{2k-1,2k}=ix_k$, when $k=1,...,n$.
 For $M\in \mathcal{H}_d$ there exists an unique $x\in \mathbb{R}^n$ such that $\{k M k^*, k\in G_d\}\cap \bar{C} =\{D(x)\}$. We call it the  radial part of $M$. These definitions of the radial part are the same as in Section \ref{Approximation}, up to some identifications.

\subsection{Gelfand Cetlin cones}
For $x, y\in \mathbb{R}^{r}$ we write $x \succeq y$ if $x$ and $y$ are interlaced, i.e.
$$x_{1}\ge y_1\ge x_2 \ge ...\ge x_r \ge y_r . $$
When $x\in \mathbb{R}^{r+1}$ we add the relation $y_r \geq x_{r+1}$. 
Let us consider the following Gelfand Cetlin cones: 
\begin{eqnarray*}GC_a&=&\{x=(x^{(1)},...,x^{(n)}): x^{(i)}\in \mathbb{R}^{i}, x^{(i)}\succeq x^{(i-1)}, i \leq n\};\\
GC_c&=&\{x=(x^{(1)},...,x^{(2n)}): x^{(2i)},x^{(2i-1)}\in \mathbb{R}_+^{i}, x^{(k)}\succeq x^{(k-1)}, i \leq n, k \leq 2n\};\\GC_b&=&\{x=(x^{(1)},...,x^{(2n)}): x^{(2i)}\in \mathbb{R}_+^{i}, x^{(2i-1)}\in \mathbb{R}_+^{i-1}\times \mathbb{R}, \vert x \vert\in GC_c\};\\ CG_d&=&\{x=(x^{(1)},...,x^{(2n-1)}): \mbox{ there exists  } x^{(2n)} \in \mathbb{R}_+^n \mbox{ such that } (x,x^{(2n)})\in CG_b\}.
\end{eqnarray*}
where, in the definition of $GC_b$, $\vert x\vert$ has the same components as $x$ except $x^{(2i-1)}_i $ which is replaced by $|x^{(2i-1)}_i|$.   The first line of the cone $GC_\nu$ is $x^{(n)}$ when $\nu=a$, $x^{(2n)}$ when $\nu=b,c$ and $x^{(2n-1)}$ when $\nu=d$. For each $\lambda$ in $ \mathbb{R}^n$ we let $GC_\nu(\lambda)$ be the cone with first line $\lambda$.\\

For $M\in \mathbb{M}_n(\mathbb{C})$ the main minor of order $m\le n$ of $M$ is the submatrix $(M_{ij})_{1\le i,j\le m}$. When $\nu=a,b,d$ the main minor of a matrix in $\mathfrak{g}_\nu$ is also in one of the classical Lie algebra (of the same type except when $\nu=b$ or $d$ and  in that case it is of type $B$ when $m$ is odd and of type $D$ when $m$ is even).  This also true for when $m$ is even when $\nu=c$. Thus we can define the radial part of a main minor  of $M\in \mathcal{H}_\nu $ as above.\\

 \begin{e-definition} For $\nu=a,b,c,d$ let $M \in \mathcal{H}_\nu$. We write 
 \begin{eqnarray*}
  \lambda_a(M)&=&(\lambda^{(1)}_a(M),...,\lambda^{(n)}_a(M)), \,\,\lambda_b(M)=(\lambda^{(2)}_b(M),\lambda^{(3)}_b(M),...,\lambda^{(2n+1)}_b(M))\\ \lambda_c(M)&=&(\lambda^{(2)}_c(M),\lambda^{(4)}_c(M),...,\lambda^{(2n)}_c(M)),\,\, \lambda_d(M)=(\lambda^{(2)}_d(M),\lambda^{(3)}_d(M),...,\lambda^{(2n)}_d(M))
 \end{eqnarray*}
where $\lambda_\nu^{(i)}(M)$ is the radial part of the main minor of order $i$ of $M$.
 \end{e-definition}$ $
 
 Using the proposition \ref{theogeneral} we get the following Theorem.\\
\begin{theorem} \label{theoclassic} For $\nu=a,b,c,d$, let  $M_\nu\in \mathcal{H}_\nu$ be a random matrix with a law invariant under the adjoint action of $G_\nu$. 
Then $\lambda_\nu(M_\nu)$, conditioned by the fact that the radial part of $M_\nu$ is $\lambda\in \mathbb{R}^n$, is uniformly distributed on $GC_\nu(\lambda)$ for $\nu=a,b,d$ and  is distributed according to the image of the uniform measure on $GC_c(\lambda)$ by the map $(x^{(1)},...,x^{(2n-1)},x^{(2n)})\in GC_c\mapsto (x^{(2)},x^{(4)},...,x^{(2n-2)},x^{(2n)})\in \mathbb{R}^{n(n+1)/2}$, for $\nu=c$.
\end{theorem}

\section{Interlaced determinantal point processes}
Considering the eigenvalues of the main minors of some random Hermitian matrices from the classical complex Lie algebras, we construct random configurations on $\mathbb{N}\times \mathbb{R}$ which verify interlacing conditions. We claim that they are determinantal and give their correlation kernels.\\
\begin{theorem} \label{theodet} For $\nu=a,b,c$ or $d$, let $M \in \mathcal{H}_\nu$ be a random matrix with a invariant law invariant under  the adjoint action of $G_\nu$. Let $\psi_{i}$, $i=1,...,n$, be measurable functions on $\mathbb{R}$, null on $\mathbb{R}_-$ for $\nu=b,c,d$, such that for all $k\in \mathbb{N},$  $x^k\psi_{i}(x)$ is integrable on $ \mathbb{R}$. We suppose that the strictly positive eigenvalues  (resp. the eigenvalues) of  $M$, $\nu=b,c,d$ (resp. $\nu=a$) have a density with respect to the Lebesgue measure proportional to 
$\Delta_\nu(x)\det(\psi_{j}(x_{i}))_{1\le i,j\le n}$, where $\Delta_a(x)=\prod_{1\le i<j\le n}(x_i-x_j)$, $\Delta_b(x)=\Delta_c(x)=\prod_{i=1}^{n}x_i\prod_{1\le i<j\le n}(x^2_i-x^2_j)$, $\Delta_d(x)=\prod_{1\le i<j\le n}(x^2_i-x^2_j)$, $x\in \mathbb{R}^n$.

Let us consider the point processes $$\xi_a=\sum_{i=1}^{n}\sum_{j=1}^{i}\delta_{i,\Lambda_{a,j}^{(i)}}, \quad \xi_b=\sum_{i=1}^{2n}\sum_{j=1}^{[(i+1)/2]}\delta_{i,\vert\Lambda_{b,j}^{(i+1)}\vert}, \quad \xi_c=\sum_{i=1}^{n}\sum_{j=1}^{i}\delta_{i,\Lambda_{c,j}^{(i)}},\quad \xi_d=\sum_{i=1}^{2n-1}\sum_{j=1}^{[(i+1)/2]}\delta_{i, \vert\Lambda_{d,j}^{(i+1)} \vert},$$
where $\Lambda_{\nu,j}^{(i)}$ is the $j^{th}$ component of $\lambda_{\nu}^{(i)}(M_\nu)$.
Then,  

$(i)$ The point process $\xi_\nu$ is  determinantal,

$(ii)$ Its  correlation kernel is
\begin{displaymath}
\begin{array}{ll} 
R((r,y),(s,z))&=-\frac{1_{s>r}(z-y\wedge z)^{\varphi_\nu(r)-\varphi_\nu(s)-1}}{(\varphi_\nu(r)-\varphi_r(s)-1)!}\\ &+\alpha_\nu\sum_{k=1}^{n}[\psi_{k}]^{-\varphi_\nu(r)}(y)\int\frac{\partial^{\varphi_\nu(s)}\Delta_ \nu}{\partial x_{k}^{\varphi_\nu(s)}}(x_1,...,x_{k-1},z,x_{k+1},...,x_n)\prod_{i\ne k}\psi_{i}(x_i)dx_i
\end{array}
\end{displaymath}
where $[\psi_{k}]^{-i}(y)=
\int_y^{+\infty}\frac{1}{(i-1)!}(x-y)^{i-1}\psi_{k}(x)\,dx$ if $i>0$, 
$[\psi_{k}]^{0}=\psi_{k}$, $\alpha_\nu^{-1}=\int \Delta_\nu(x)\prod_{i=1}^{n}\psi_{i}(x_i)dx_i$, and $\varphi_a(r)=\frac{1}{2}\varphi_c(r)=n-r$, $\varphi_b(r)=\varphi_d(r)+1=2n-r$,

$(iii)$ Moreover, when we can write $\Delta_\nu(x)=\det(\chi_{i}(x_{j}))_{1\le i,j\le n}$ where $(\chi_{i})_{i\ge 1}$ is a sequence of functions on $\mathbb{R}$ such that $\chi_{i}\psi_{j}$ is integrable on $\mathbb{R}$ and $\int_\mathbb{R} \chi_{i}(x)\psi_{j}(x)dx=\delta_{ij}$, $i,j=1,...,n$, then
$$
R_\nu((r,y),(s,z))=-\frac{1_{s>r}(z-y\wedge z)^{\varphi_\nu(r)-\varphi_\nu(s)-1}}{(\varphi_\nu(r)-\varphi_\nu(s)-1)!}+\sum_{k=1}^{n}[\psi_{k}]^{-\varphi_\nu(r)}(y)\frac{d^{\varphi_\nu(s)}\chi_{k}}{d x^{\varphi_\nu(s)}}(z) .
$$
\end{theorem}

Let us describe some applications of this theorem:
\medskip

\noindent {\bf The Gaussian case} For $\nu=a,b,c,d$, let $M_\nu \in \mathcal{H}_\nu$ be a random Gaussian matrix distributed according to a probability measure proportional to $e^{-t_\nu Tr(H^2)}\mu_\nu(dH)$ where $\mu_\nu$ is the Lebesgue measure on $\mathcal{H}_\nu$, $t_\nu=\frac{1}{2}$ for $\nu=b,c,d$ and $t_a=1$.
The matrices $M_\nu$, $\nu=a,d,b,c$, satisfy the hypothesis of the Theorem by taking respectively $\nu=a$ and $\psi_i(x)=x^{i-1}e^{-x^2},$, $\nu=d$ and $\psi_i(x)=x^{2i-2}e^{-x^2}1_{x>0}$, and $\nu=b,c$ and  $\psi_i^\nu(x)=x^{2i-1}e^{-x^2}1_{x>0},$. Besides, the hypothesis of the point $(iii)$ are satisfied if we chose $\chi_i=h_{i-1}$ for $\nu=a$, $\chi_i=h_{2i-1}$, for $\nu=b,c$,  and $\chi_i=h_{2i-2}$, for $\nu=d$, where $(h_i)_{i\ge 0}$ is the sequence of Hermite normalised  polynomials such that $h_i$ has degree $i$.

\medskip
\noindent{\bf GUE and LUE Ensembles} The Gaussian, Laguerre and Jacobi unitary ensembles  are obtained by taking $\nu=a$, $\psi_i(x)=x^{i-1} e^{-\alpha x^2},$ $\psi_i(x)=x^{i-1}x^\alpha e^{-\beta x}1_{x>0},$ and  $\psi_i(x)=x^{i-1} x^{\alpha}(1-x)^{\beta}1_{0<x<1}$. 
\medskip

 If the radial part of $M$ is  deterministic and equal to $\lambda\in \mathbb{R}^n$, the theorem  remains true up to slight modifications, replacing $\psi_{i}(x)dx$ by $\delta_{\vert \lambda_i\vert}(dx)$ in the kernel $R$. As we have seen in the example above,  $(iii)$ of  the Theorem \ref{theodet} generalises Theorem $1.3$ of \cite{Johansson}. Let us give a similar result  for the orthogonal case.

\begin{corollary} Let $M$ be distributed according to a probability measure proportional to  $e^{-\frac{1}{2}Tr(H^2)}dH$, where $dH$ is the Lebesgue measure on the anti-symmetric Hermitian matrices. Consider the random vectors $\Lambda^{(i)}\in \mathbb{R}^{[(i+1)/2]}$ of strictly positive eigenvalues of its main minor of order $i+1$, $i\in \mathbb{N}^*$. Then the point process  $\sum_{i=1}^{+\infty}\sum_{j=1}^{[(i+1)/2]}\delta_{i, \Lambda^{(i)}_{j}}$ is determinantal on $\mathbb{N}^*\times \mathbb{R}_{+}$ with correlation kernel
\begin{displaymath}
\begin{array}{ll} 
R((r,y),(s,z))=&-\frac{1_{r<s}}{(s-r-1)!} (z-y)^{s-r-1}1_{y<z}\\&+\sum_{i=0}^{[\frac{r+1}{2}]\wedge [\frac{s+1}{2}]}\frac{(2^{r}(r-2i+1)!)^{1/2}}{(2^{s}(s-2i+1)!)^{1/2}}h_{s-2i+1}(z)h_{r-2i+1}(y)e^{-y^{2}} \\
&+ \sum_{i=[\frac{r+1}{2}]+1}^{[\frac{s+1}{2}]}\frac{h_{s-2i+1}(z)}{(2^{s-2i+1}(s-2i+1)!\sqrt{\pi})^{1/2}}\int^{+\infty}_{y} \frac{(x-y)^{2i-r-2}}{(2i-r-2)!}e^{-x^{2}}dx
\end{array}
\end{displaymath}
\end{corollary}

\paragraph*{}For types $A$ and $C$, the proof of Theorem \ref{theodet} rests  on the criterion of lemma 3.4 in  \cite{Borodin}, adapted for a continuous framework. For the orthogonal cases  $B$ and $D$, the criterion needs to be slightly modified. For instance, for the odd orthogonal case, we show that:
\\

\begin{e-proposition}   Let $\phi_{2r-1}:\mathbb{R}_+\times \mathbb{R}_+ \to \mathbb{R}_+$, $\phi_{2r-2}:\mathbb{R}_+\times \mathbb{R}_+ \to \mathbb{R}_+$, and $\psi_{r}:\mathbb{R}_+ \to \mathbb{R}_+$, $r=1,...,n$,  be measurable functions.   Let $X^{(2i-1)}\in \mathbb{R}_+^i$ and $X^{(2i)}\in \mathbb{R}_+^i$, $i=1,...,n$, be $2n$ random variables. Suppose that $(X^{(1)},...,X^{(2n)})$  has a density with respect to the Lebesgue measure proportional to
\begin{displaymath}
\begin{array}{l} 
\prod_{r=1}^{n}[\det(\phi_{2r-2}(x_i^{(2r-2)},x_j^{(2r-1)}))_{i,j\le r}\det(\phi_{2r-1}(x_i^{(2r-1)},x_j^{(2r)}))_{i,j\le r}]  \det(\psi_{j}(x_i^{(2n)}))_{1\le i,j\le n},
\end{array}
\end{displaymath}
where $x_r^{2r-2}\in \mathbb{R}$, $r=1,...,n$ are arbitrary fixed real numbers.
 Let us write 
\begin{displaymath}
\begin{array}{l} 
\phi_{r}*\phi_{s}(x,y)=\int_{0}^\infty \phi_{r}(x,z)\phi_{s}(z,y)dz
\end{array},
\quad \phi^{(r,s)}=\left\{\begin{array}{ll} 
\phi_{r}*\phi_{r+1}*...*\phi_{s-1}& \textrm{ if $r<s$}\\
0 & \textrm{ otherwise.}
\end{array}\right.
\end{displaymath}

\begin{displaymath}
\phi_{r-1}*\phi^{(r,s)}=\left\{\begin{array}{ll} 
 \phi_{r-1}*\phi^{(r,s)}& \textrm{if $r<s$}\\
\phi_{r-1}& \textrm{if $r=s$}
\end{array}\right., \quad
\psi_{r-k}^{r}(x)=\left\{\begin{array}{ll} 
\int_{0}^{\infty}  \phi^{(r,2n)}(x,y)\psi_k(y)dy & \textrm{if $r<2n$.}\\
\psi_k(x)& \textrm{if $r=2n$}
\end{array}\right.
\end{displaymath}
Let $M=(M_{ij})_{1\le i,j\le n}$ defined by 
$M_{ij}=\int_{0}^{\infty} \phi_{2i-2}*\phi^{(2i-1,2n)}(x_{i}^{2i-2},x)\psi_{j}(x)dx$. Suppose that all these integrals are finite. 
Then $M$ is invertible and the point process $\sum_{i=1}^{2n}\sum_{j=1}^{[(i+1)/2]}\delta_{i,X_j^{(i)}}$ is  determinantal on $\{1,...,2n\}\times \mathbb{R}_+$ with correlation kernel 
\begin{displaymath}
\begin{array}{l} 
K((r,x),(s,y))=-\phi^{(r,s)}(x,y)+\sum_{k=1}^{n} \psi_{r-k}^{r}(x)\sum_{l=1}^{[(s+1)/2]}(M^{-1})_{kl}\phi_{2l-2}*\phi^{(2l-1,s)}(x_l^{2l-2},y).
\end{array}
\end{displaymath}
\end{e-proposition}
\begin{remark}
The point process $\xi_b$ defined at the Theorem \ref{theodet} is obtained by taking  $x_r^{2r-2}=0$,  $r=1,...,n$, and $\phi_{r}(x,y)=1_{y\ge x}$, for $x,y\in \mathbb{R}_+$, $r=1,...,2n$. 
\end{remark}



\end{document}